\documentclass[a4paper,12pt]{article}

\usepackage{amsmath,amssymb,amsthm}
\usepackage{mathtools}
\usepackage[dvips,final]{graphicx}
\usepackage[dvips,final]{epsfig}

\usepackage{hyperref}
\hypersetup{
	colorlinks=true,
	linkcolor=blue,
	urlcolor=cyan,
	citecolor=blue,
}

\usepackage{xcolor}

\definecolor{oran}{rgb}{.75,.4,0}
\definecolor{vio}{rgb}{0.7,0,0.5}
\definecolor{gre}{rgb}{0.1,0.6,0}
\definecolor{ora}{rgb}{0.8,0.2,0}
\definecolor{blgr}{rgb}{0,0.5,0.5}

\usepackage{titlesec}
\titleformat{\section}{\bfseries}{\thesection}{1em}{}
\titleformat{\subsection}{\itshape}{\thesubsection}{1em}{}

\numberwithin{equation}{section}

\newfont{\ctv}{msam10}

\newcommand{\bbox}{\mbox{\ctv \symbol{4}}}
\def\QED{{${}\hfill\bbox$}}
\newenvironment{pf}[1]{\par\vskip1mm{\noindent\it #1.}\ }{\QED\par
\vskip2mm}

\def\bpf{\begin{pf}}
\def\epf{\end{pf}}

\def\expe{\hbox{\rm e}}
\def\ve{\varepsilon}
\def\vp{\varphi}
\def\dd{\,\mathrm{d}}

\def\sign{\mathrm{\,sign}}
\def\supess{\mathop{\mathrm{\,sup\,ess}}}
\def\for{\mathrm{\ for\ }}
\def\ale{\mathrm{\ a.\,e.}}

\def\play{\mathfrak{p}}
\def\on{^{(n)}}
\def\ota{^{(\tau)}}

\def\sumi{\sum_{i=0}^{\infty}}

\def\real{\mathbb{R}}
\def\nat{\mathbb{N}}

\def\io{\int_{\Omega}}

\def\be{\begin{equation}\label}
\def\ee{\end{equation}}

\def\ber{\begin{eqnarray}}
\def\eer{\end{eqnarray}}

\def\bers{\begin{eqnarray*}}
\def\eers{\end{eqnarray*}}

\def\bpf{\begin{pf}}
\def\epf{\end{pf}}

\newtheorem{theorem}{Theorem}[section]

\newtheorem{hypothesis}[theorem]{Hypothesis}
\newtheorem{proposition}[theorem]{Proposition}
\newtheorem{definition}[theorem]{Definition}

\begin{document}

\title{Long-time behaviour of a porous medium model with degenerate hysteresis
\thanks{The support from the Austrian Science Fund (FWF) grants V662, Y1292, and F65, from the OeAD WTZ grants CZ02/2022, CZ09/2023, and MULT06/2023, and from the M\v{S}MT grants 8J22AT017, 8J23AT008, and 8X23001 is gratefully acknowledged.}
}

\author{Chiara Gavioli
\thanks{Institute of Analysis and Scientific Computing, TU Wien, Wiedner Hauptstra\ss e 8-10, A-1040 Vienna, Austria, E-mail: {\tt chiara.gavioli@tuwien.ac.at}.}
\and Pavel Krej\v c\'{\i}
\thanks{Faculty of Civil Engineering, Czech Technical University, Th\'akurova 7, CZ-16629 Praha 6, Czech Republic, E-mail: {\tt Pavel.Krejci@cvut.cz}.}
}

\date{}

\maketitle

\begin{abstract}
Hysteresis in the pressure-saturation relation in unsaturated porous media, owing to surface tension on the liquid-gas interface, exhibits strong degeneracy in the resulting mass balance equation. As an extension of previous existence and uniqueness results, we prove that under physically admissible initial conditions and without mass exchange with the exterior, the unique global solution of the fluid diffusion problem exists and asymptotically converges as time tends to infinity to a possibly non-homogeneous mass distribution and an a priori unknown constant pressure.

\bigskip


\medskip

\end{abstract}


\section{Introduction}\label{sec:intr}

This article deals with the problem of existence, uniqueness, and long-time stabilization of the solution to the degenerate PDE with hysteresis in a bounded space domain $\Omega \subset \real^N$, $N \in \nat$ of class $C^{1,1}$
\be{e0}
s_t - \Delta u = 0 \quad \for (x,t) \in \Omega\times (0,\infty)
\ee
where $t>0$ is the time variable, $u = u(x,t) \in \real$ represents the pressure, $s = s(x,t) \in (0,1)$ is the relative saturation of the fluid in the pores, and $\Delta$ is the Laplacian in $x$. Hysteresis in the pressure-saturation relation is represented by a Preisach operator $G$ in the form
\be{e0a}
s(x,t) = G[u](x,t).
\ee
A detailed justification of why it is meaningful to consider Preisach hysteresis in porous media modelling can be found in \cite{vis} or in the introduction of \cite{colli}.

System equations \eqref{e0} and \eqref{e0a} are considered with the boundary condition
\be{e1}
-\nabla u(x,t)\cdot n = 0
\ee
on $\partial\Omega$, and with a given initial condition that includes not only the initial pressure
\be{e2}
u(x,0) = u_0(x),
\ee
but also an initial Preisach memory distribution specified below in a rigorous formulation of the Preisach operator in Definition~\ref{dpr}. Let us only mention at this point that the time evolution described by equation~\eqref{e0} is doubly degenerate: in typical situations, the function $s(x,t) = G[u](x,t)$ is bounded independently of the evolution of $u$, so that no a priori lower bound for $u_t$ is immediately available. Furthermore, at every point $x \in \Omega$ and every time $t_0$ where $u_t$ changes sign (which the engineers call a {\em turning point\/}), we have
\be{e6a}
u_t(x,t_0-\delta)\cdot u_t(x,t_0+\delta)<0 \ \ \forall \delta\in (0,\delta_0(x)) \ \Longrightarrow \ \liminf_{\delta \to 0+} \frac{G[u]_t(x,t_0+\delta)}{u_t(x,t_0+\delta)} = 0,
\ee
so that the knowledge of $G[u]_t$ alone does not give a complete information about $u_t$. In particular, if $\Delta u_0(x) \ne 0$ and the initial memory of $G$ corresponds to a turning point, then even a local solution cannot be expected to exist. Hypothesis~\ref{hym} is shown to avoid this pathological situation. A more detailed discussion about this issue can be found in the recent paper \cite{colli}. There, under suitable hypotheses on the data, the existence and uniqueness of a solution to equation~\eqref{e0} on an arbitrary time interval $t\in (0,T)$ has been proved in the case of Robin boundary conditions. Crucial assumptions to obtain the existence result are that the operator $G$ is a so-called convexifiable Preisach operator, and its initial memory state is compatible with the initial condition $u_0$ in equation~\eqref{e2} (see Section~\ref{sec:stat}). The present autonomous case with homogeneous Neumann boundary conditions makes it possible to prove a stronger existence and uniqueness statement (see the first paragraph after Theorem~\ref{t2}) under the same hypotheses on the data, as well as asymptotic convergence to an a priori unknown equilibrium as time tends to infinity.

The structure of the paper is as follows. In Section~\ref{sec:stat}, we list the definitions of the main concepts, including convexifiable Preisach operators, and state the main Theorem~\ref{t2}. In Section~\ref{disc}, we propose a time discretization scheme with time step $\tau>0$, and in Sections~\ref{unif}--\ref{conv}, we derive estimates independent of $\tau$. In particular, estimate \eqref{fest} improves the corresponding convexity estimate obtained in \cite{colli}. By a compactness argument, we pass to the limit as $\tau \to 0$ and prove in Section~\ref{limi} that the limit is the unique solution to our PDE problem with hysteresis. The long-time stabilization result is proved in Section~\ref{sec:asym}.


\section{Statement of the problem}\label{sec:stat}

The Preisach operator was originally introduced in \cite{prei}. For our purposes, it is more convenient to use the equivalent variational setting from \cite{book}.

\begin{definition}\label{dpr}
	Let $\lambda \in L^\infty(\Omega \times (0,\infty))$ be a given function which we call the {\em initial memory distribution\/} and which has the following properties: 
	\begin{align}\label{e6}
		&|\lambda(x, r_1) - \lambda(x, r_2)| \le |r_1 - r_2| \ \ale \ \forall r_1, r_2 \in (0,\infty),\\ \label{e6b}
		&\exists \Lambda > 0:\ \lambda(x,r) = 0\ \for r\ge \Lambda \ \mbox{ and a.\,e.\ } x \in \Omega.
	\end{align}
	For a given $r>0$, we call the {\em play operator with threshold $r$ and initial memory $\lambda$} the mapping which, with a given function $u \in L^p(\Omega; W^{1,1}(0,T))$ for $p \ge 1$, associates the solution $\xi^r\in L^p(\Omega; W^{1,1}(0,T))$ of the variational inequality
	\be{e4a}
	|u(x,t) - \xi^r(x,t)| \le r, \quad \xi^r_t(x,t)(u(x,t) - \xi^r(x,t) - z) \ge 0 \ \ale \ \forall z \in [-r,r],
	\ee
	with a given initial memory distribution
	\be{e5}
	\xi^r(x,0) = \lambda(x,r),
	\ee
	and we denote
	\be{e4}
	\xi^r(x,t) = \play_r[\lambda,u](x,t).
	\ee
	Given a measurable function $\rho :\Omega\times(0,\infty)\times \real \to [0,\infty)$ and a constant $\bar G \in \real$, the Preisach operator $G$ is defined as a mapping $G: L^p(\Omega; W^{1,1}(0,T))\to L^p(\Omega; W^{1,1}(0,T))$ by the formula
	\be{e3}
	G[u](x,t) = \bar G + \ \int_0^\infty\int_0^{\xi^r(x,t)} \rho(x,r,v)\dd v\dd r.
	\ee
	The Preisach operator is said to be {\em regular\/} if the density function $\rho$ of $G$ in \eqref{e3} belongs to $L^1(\Omega\times (0,\infty )\times \real)\cap L^\infty(\Omega\times (0,\infty )\times \real)$, and there exist a constant $\rho_1$ and a decreasing function $\rho_0: \real\to \real$ such that for all $U>0$ we have
	\be{e3a}
	\rho_1 > \rho(x,r,v)> \rho_0(U) > 0 \ \for \ale\ (x,r,v)\in \Omega\times (0,U) \times (-U,U).
	\ee
\end{definition}

Let us mention the following classical result (see \cite[Proposition~II.3.11]{book}).

\begin{proposition}\label{pc1}
	Let $G$ be a regular Preisach operator in the sense of Definition~\ref{dpr}. Then it can be extended to a Lipschitz continuous mapping $G: L^p(\Omega; C[0,T]) \to L^p(\Omega; C[0,T])$ for every $p \in [1,\infty)$.
\end{proposition}

The Preisach operator is rate-independent. Hence, for input functions $u(x,t)$ which are monotone in a time interval $t\in (a,b)$, a regular Preisach operator $G$ can be represented by a superposition operator $G[u](x,t) = B(x, u(x,t))$ with an increasing function $u \mapsto B(x, u)$ called a {\em Preisach branch\/}. Indeed, the branches may be different at different points $x$ and different intervals $(a,b)$. The branches corresponding to increasing inputs are said to be {\em ascending\/} (the so-called wetting curves in the context of porous media), the branches corresponding to decreasing inputs are said to be {\em descending\/} (drying curves).

\begin{definition}\label{dpc}
	Let $U>0$ be given. A Preisach operator is said to be {\em uniformly counterclockwise convex on $[-U,U]$\/} if for all inputs $u$ such that $|u(x,t)|\le U$ a.\,e., all ascending branches are uniformly convex and all descending branches are uniformly concave.
	
	A regular Preisach operator is called {\em convexifiable\/} if for every $U>0$ there exist a uniformly counterclockwise convex Preisach operator $P$ on $[-U,U]$, positive constants $g_*(U),g^*(U),\bar{g}(U)$, and a twice continuously differentiable mapping $g:[-U,U] \to [-U,U]$ such that
	\begin{equation}\label{g}
		g(0)=0, \quad 0<g_*(U) \le g'(u) \le g^*(U), \quad |g''(u)| \le \bar g(U)\ \ \forall u\in [-U,U],
	\end{equation}
	and $G = P\circ g$.
\end{definition}

A typical example of a uniformly counterclockwise convex operator is the so-called {\em Prandtl-Ishlinskii operator\/} characterized by positive density functions $\rho(x,r)$ independent of $v$, see \cite[Section~4.2]{book}. Operators of the form $P\circ g$ with a Prandtl-Ishlinskii operator $P$ and an increasing function $g$ are often used in control engineering because of their explicit inversion formulas, see \cite{al,viso,kk}. They are called the {\em generalized Prandtl-Ishlinskii operators\/} (GPI) and represent an important subclass of Preisach operators. Note also that for every Preisach operator $P$ and every Lipschitz continuous increasing function $g$, the superposition operator $G = P\circ g$ is also a Preisach operator, and there exists an explicit formula for its density, see \cite[Proposition~2.3]{error}. Another class of convexifiable Preisach operators is shown in \cite[Proposition~1.3]{colli}. 

As it has been mentioned in Section~\ref{sec:intr}, even a local solution to Problem~\eqref{e0}--\eqref{e2} may fail to exist if for example $\lambda(x,r) \equiv 0$ and $\Delta u_0(x)\ne 0$. Then $t=0$ is a turning point for all $x \in \Omega$ and there is no way to satisfy equation~\eqref{e0} in any sense. We therefore need an initial memory compatibility condition which we state in the following way.

\begin{hypothesis}\label{hym}
	The initial pressure $u_0$ belongs to $W^{2,2}(\Omega)$, $\Delta u_0 \in L^\infty(\Omega)$, and there exist constants $L>0$, $\Lambda > 0$ and a function $r_0 \in L^\infty(\Omega)$ such that $\supess_{x\in \Omega}|u_0(x)| \le \Lambda$, equation~\eqref{e6b} is satisfied, and the following initial compatibility conditions hold:
	\begin{align} \label{c0}
		\lambda(x,0) &= u_0(x) \ \ale \textup{ in } \Omega,\\ \label{c1}
		\frac1L \sqrt{|\Delta u_0(x)|} &\le r_0(x) \le \Lambda \ \ale \textup{ in } \Omega,
		\\ \label{c2} -\frac{\partial}{\partial r} \lambda(x,r) &\in \sign(\Delta u_0(x)) \ \ale \textup{ in } \Omega\ \for r\in (0,r_0(x)), \\ \label{c2a}
		-\nabla u_0(x)\cdot n &= 0 \ \ale\ \mbox{\em on }\, \partial\Omega.
	\end{align}
\end{hypothesis}

It was shown in \cite[Proposition~2.2.16]{BS} that the solution to inequality~\eqref{e4a} is such that for every $u\in L^p(\Omega;W^{1,1}(0,T))$ the property
\begin{equation}\label{lip-r}
	|\xi^{r_1}(x, t) - \xi^{r_2}(x,t)| \le |r_1 - r_2| \ \ale \ \ \forall r_1, r_2 \in (0,\infty)
\end{equation}
is preserved during the evolution. From relations \eqref{e6}--\eqref{e5} it follows that for a.\,e.\ $(x,t) \in \Omega\times (0,\infty)$ there exists a so-called {\em active memory level\/} $r^*(x,t) \ge 0$ such that $|\xi^r(x,t) - u(x,t)| = r$ and $\xi^r_t(x,t) = u_t(x,t)$ for $r<r^*(x,t)$, $|\xi^r(x,t) - u(x,t)| < r$ and $\xi^r_t(x,t) = 0$ for $r>r^*(x,t)$. In this sense, the value $r_0$ in conditions \eqref{c1} and \eqref{c2} represents the active memory level at time $t=0$. The meaning of conditions \eqref{c0}--\eqref{c2} is that for large values of $\Delta u_0(x)$, the initial memory $\lambda$ has to go deeper in the memory direction. We refer to \cite{colli} for explanation how Hypothesis~\ref{hym} guarantees the existence of some previous admissible history of the process prior to the time $t=0$ and ensures the existence of a continuation for $t>0$.

Problem \eqref{e0}--\eqref{e1} is to be understood in variational form
\be{e0v}
\io G[u]_t\vp \dd x + \io\nabla u\cdot\nabla \vp\dd x = 0
\ee
for every test function $\vp \in W^{1,2}(\Omega)$. Existence and uniqueness of a solution to equations \eqref{e0} and \eqref{e0a} with the Robin boundary condition and admissible initial conditions has been proved in \cite{colli}. We now state the main result of this paper.

\begin{theorem}\label{t2}
	Let Hypothesis~\ref{hym} be satisfied, and let $G$ be a convexifiable Preisach operator in the sense of Definition~\ref{dpc}. Then Problem~\eqref{e0v} with initial condition $u(x,0) = u_0(x)$, $\play_r[\lambda,u](x,0) = \lambda(x,r)$ for all $r>0$ admits a unique global solution $u$ on $\Omega\times (0,\infty)$ such that
		$$
		\begin{aligned}
			u &\in L^\infty(\Omega\times (0,\infty)), \quad |u(x,t)| \le \Lambda \ \ale,\\[-.25em]
			\nabla u &\in L^\infty(0,\infty; L^2(\Omega; \real^N)) \cap L^2(\Omega\times (0,\infty);\real^N),\\[-.25em]
			\Delta u, G[u]_t &\in L^2(\Omega\times (0,\infty)),\\[-.25em]
			u_t &\in L^Q(\Omega\times (0,\infty)),
		\end{aligned}
		$$
		where $L^Q(\Omega\times (0,\infty))$ is the Orlicz space generated by a function $Q(z)$ which behaves like $z^3$ for $z\in (0,1)$ and $z^p$ with $p = 1+ (2/N)$ for $z > 1$, see equation~\eqref{defq} for $N\ge 3$, equation~\eqref{Q2} for $N=2$, and equation~\eqref{Q1} for $N=1$. Moreover, there exist a constant $\bar u \in \real$ such that $\lim_{t\to \infty}\io |u(x,t) - \bar u|^q \dd x = 0$ for all $q\ge 1$, and a function $\bar\lambda\in L^\infty(\Omega\times [0,\infty))$ such that $\bar\lambda(x,0) = \bar u$, $|\bar\lambda(x,r)-\bar\lambda(x,s)| \le r-s$ for a.\,e.\ $x\in\Omega$ and all $r>s\ge 0$, and such that $\lim_{t\to \infty}\io |\play_r[\lambda,u](x,t) - \bar \lambda(x,r)|^q \dd x = 0$ for all $q\ge 1$ and $r\ge 0$.
\end{theorem}

Note that the regularity of $u_t$ in \cite{colli} is only $L^{p-\ve}$ for each $\ve > 0$ and on bounded time intervals. We propose here a refined estimation technique which allows to derive a global in time Orlicz bound for $u_t$(for the theory of Orlicz spaces, refer to \cite{kala}) and full $L^p$ bound on bounded time intervals.

Putting $\vp = 1$ in equation~\eqref{e0v} we formally get the identity
$$
\frac{\dd}{\dd t}\io s(x,t)\dd x = 0,
$$
which means that the total liquid mass is preserved during the evolution. The meaning of Theorem~\ref{t2} is that the pressure $u$ is asymptotically uniformly distributed in $\Omega$, but because of hysteresis, we cannot expect that the spatial distribution of the liquid mass will also be asymptotically uniform and that the constant limit pressure value $\bar u$ can be computed explicitly.


\section{Time discretization}\label{disc}

We proceed as in \cite{colli}, choose a sufficiently small time step $\tau > 0$, and replace equation~\eqref{e0v} with its time discrete system for the unknowns $\{u_i: i \in \nat\cup\{0\}\} \subset W^{1,2}(\Omega)$ of the form
\be{dis1}
\io \left(\frac1\tau(G[u]_i - G[u]_{i-1})\vp + \nabla u_i\cdot\nabla\vp\right)\dd x = 0
\ee
for $i\in \nat$ and for every test function $\vp \in W^{1,2}(\Omega)$, where $u_0$ is the initial condition in equation~\eqref{e2}. Here, the time-discrete Preisach operator $G[u]_i$ is defined by a formula of the form \eqref{e3}, namely,
\be{de3}
G[u]_i(x) = \bar G + \ \int_0^{\infty}\int_0^{\xi^r_i(x)} \rho(x,r,v)\dd v\dd r,
\ee
where $\xi^r_i$ denotes the output of the time-discrete play operator
\be{de4}
\xi^r_i(x) = \play_r[\lambda,u]_i(x)
\ee
defined as the solution operator of the variational inequality
\be{de4a}
|u_i(x) - \xi^r_i(x)| \le r, \quad (\xi^r_i(x) - \xi^r_{i-1}(x))(u_i(x) - \xi^r_i(x) - z) \ge 0 \quad \forall i\in \nat \ \ \forall z \in [-r,r],
\ee
with a given initial memory curve
\be{de5}
\xi^r_0(x) = \lambda(x,r) \ \ale,
\ee
in a similar manner to relations \eqref{e4a} and \eqref{e5}. Note that the discrete variational inequality \eqref{de4a} can be interpreted as weak formulation of \eqref{e4a} for piecewise constant inputs in terms of the Kurzweil integral, and details can be found in \cite[Section 2]{ele}.

Arguing as in \cite[Eq.~(35)]{colli}, we obtain the two-sided estimate
\be{de5a}
\frac1C|G[u]_i(x) - G[u]_{i-1}(x)|^2 \le (u_i(x) - u_{i-1}(x))(G[u]_i(x) - G[u]_{i-1}(x)) \le C|u_i(x) - u_{i-1}(x)|^2,
\ee
with a constant $C>0$ depending only on the constant $\rho_1$ from relation~\eqref{e3a}.

For each $i\in\nat$, there is no hysteresis in the passage from $u_{i-1}$ to $u_i$, so that equation~\eqref{dis1} is a standard 
monotone semilinear elliptic equation which admits a unique solution $u_i \in W^{1,2}(\Omega)$ for every $i \in \nat \cup \{0\}$.


\section{Uniform upper bound}\label{unif}

Following \cite{hilp}, the idea is to test \eqref{dis1} by $\vp =H_\ve(u_i - \Lambda)$, with $H_\ve$ being a Lipschitz regularization of the Heaviside function $H(s) = 1$ for $s>0$, $H(s) = 0$ for $s\le 0$, and $\Lambda$ from \eqref{e6b}. We then let $\ve$ tend to $0$. The elliptic term gives a non-negative contribution, and we get for all $i \in \nat$ that
\be{hil}
\io (G[u]_i - G[u]_{i-1})H(u_i - \Lambda)\dd x \le 0.
\ee
We define the functions
\be{psi}
\psi(x,r,\xi) \coloneqq \int_0^\xi\rho(x,r,v)\dd v, \quad \Psi (x,r,\xi) \coloneqq \int_0^\xi \psi(x,r,v)\dd v.
\ee
In terms of the sequence $\xi^r_i(x) = \play_r[\lambda,u]_i$ we have
\begin{equation}\label{Gi}
	G[u]_i(x) = \bar G + \int_0^\infty \psi(x,r,\xi^r_i(x))\dd r.
\end{equation}
Choosing $z= \Lambda - (\Lambda-r)^+ = \min\{\Lambda,r\}$ in \eqref{de4a} and using the fact that $\psi$ is an increasing function of $\xi$, we get for all $i\in\nat$, all $r>0$, and a.\,e.\ $x \in \Omega$ that
$$
\big(\psi(x,r,\xi^r_i(x)) - \psi(x,r,\xi^r_{i-1}(x))\big)\big((u_i(x) - \Lambda) - (\xi^r_i(x) - (\Lambda-r)^+)\big) \ge 0.
$$
The Heaviside function $H$ is non-decreasing, hence,
$$
\big(\psi(x,r,\xi^r_i(x)) - \psi(x,r,\xi^r_{i-1}(x))\big)\big(H(u_i(x) - \Lambda) - H(\xi^r_i(x) - (\Lambda-r)^+)\big) \ge 0.
$$
From \eqref{hil} it follows that
\be{hil2}
\io\int_0^\infty\big(\psi(x,r,\xi^r_i(x)) - \psi(x,r,\xi^r_{i-1}(x))\big) H(\xi^r_i(x) - (\Lambda-r)^+)\dd r \dd x \le 0.
\ee
We have by relations \eqref{e6} and \eqref{e6b} that $\xi^r_0(x) = \lambda(x,r) \le (\Lambda - r)^+$ a.\,e. We now proceed by induction assuming that $\xi^r_{i-1}(x) \le (\Lambda-r)^+$ a.\,e. for some $i\in \nat$. By \eqref{hil2} we have
\be{hil3}
\io\int_0^\infty\big(\psi(x,r,\xi^r_i(x)) {-} \psi(x,r,\xi^r_{i-1}(x))\big) \big(H(\xi^r_i(x) {-} (\Lambda{-}r)^+) {-} H(\xi^r_{i-1}(x) - (\Lambda{-}r)^+)\big)\dd r \dd x \le 0.
\ee
The expression under the integral sign in \eqref{hil3} is non-negative almost everywhere, hence it vanishes almost everywhere, and we conclude that $\xi^r_i(x) \le (\Lambda - r)^+$ a.\,e.\ for all $r\ge 0$ and $i\in \nat$. Similarly, putting $z= (\Lambda-r)^+ - \Lambda = -\min\{\Lambda,r\}$ in \eqref{de4a} we get $\xi^r_i(x) \ge -(\Lambda - r)^+$, so that
\be{Lam}
|u_i(x)| \le \Lambda, \quad |\xi^r_i(x)| \le (\Lambda - r)^+
\ee
for a.\,e.\ $x \in \Omega$ and all $r\ge 0$ and $i\in \nat$. In particular, this implies that in \eqref{Gi} we can actually write
\begin{equation}\label{Gi_Lambda}
	G[u]_i(x) = \bar G + \int_0^\Lambda \psi(x,r,\xi^r_i(x))\dd r.
\end{equation}


\section{Estimates of the pressure}\label{pres}

We first test \eqref{dis1} by $\vp = u_i - u_{i-1}$ and get
\be{des2}
\frac1{\tau}\io (G[u]_i - G[u]_{i-1})(u_i - u_{i-1})\dd x + \io \nabla u_i\cdot\nabla(u_i - u_{i-1})\dd x = 0
\ee
for all $i\in\nat$. Using the elementary inequality $\nabla u_i\cdot\nabla(u_i - u_{i-1}) \ge \frac12\left(|\nabla u_i|^2 - |\nabla u_{i-1}|^2\right)$ and putting $V_i = \frac12 \io |\nabla u_i|^2\dd x$ we obtain
\be{des3}
\frac1{\tau}\io (G[u]_i - G[u]_{i-1})(u_i - u_{i-1})\dd x +V_i - V_{i-1} \le 0.
\ee
By comparison in \eqref{dis1} (that is, testing by a function $\varphi \in W^{1,2}(\Omega)$ with compact support in $\Omega$ and using the fact that such functions form a dense subset of $L^2(\Omega)$), and by employing estimate \eqref{de5a}, we get for all $i\in\nat$ that
\be{ea2}
\io|\Delta u_i|^2 \dd x = \frac1{\tau^2}\io |G[u]_i - G[u]_{i-1}|^2 \dd x
\le \frac{c}{\tau^2}\io (G[u]_i - G[u]_{i-1})(u_i - u_{i-1})\dd x
\ee
with a constant $c>0$ independent of $i$. Coming back to \eqref{des3}, for all $i\in\nat$ we thus have
\be{ea2a}
\io|\Delta u_i|^2 \dd x + \frac{c}{\tau} (V_i - V_{i-1}) \le 0.
\ee
Consider the complete orthonormal basis $\{e_k: k\in \nat\}$ in $L^2(\Omega)$ of eigenfunctions of the operator
\be{ek}
-\Delta e_k = \mu_k e_k\ \ \mbox{ in }\ \Omega, \quad -\nabla e_k\cdot n = 0 \ \mbox{ on }\ \partial\Omega
\ee
with eigenvalues $0 = \mu_0 < \mu _1 \le \mu_2 \le \dots$. The Fourier expansion of $u_i$ in terms of the basis $\{e_k\}$ is of the form
\be{ea3}
u_i(x) = \sum_{k=1}^\infty u_k^i e_k(x)
\ee
for all $i \in \nat \cup \{0\}$. Using the orthonormality of the system $\{e_k\}$, we rewrite \eqref{ea2a} in the form
\be{ea4}
\sum_{k=1}^\infty \mu_k^2 \left(u_k^i\right)^2 + \frac{c}{\tau}\sum_{k=1}^\infty \mu_k \left((u_k^i)^2 - (u_k^{i-1})^2\right) \le 0.
\ee
The sequence
$$
V_i = \io|\nabla u_i|^2\dd x = \sum_{k=1}^\infty \mu_k \left(u_k^i\right)^2
$$
thus satisfies the inequality
$$
\mu_1 V_i + \frac{c}{\tau} (V_i - V_{i-1}) \le 0,
$$
and we conclude that
\be{vi}
V_i \le V_0 \left(1+ \frac{\mu_1\tau}{c}\right)^{-i},
\ee
that is, $V_i$ decays exponentially as $i \to \infty$.


\section{Convexity estimate}\label{conv}

By \eqref{Lam}, the functions $u_i$ do not leave the interval $[-\Lambda,\Lambda]$. Since $G$ is convexifiable in the sense of Definition~\ref{dpc}, there exist positive numbers $g_*,g^*,\bar{g}$ and a twice continuously differentiable mapping $g:[-\Lambda,\Lambda] \to [-\Lambda,\Lambda]$ such that $g(0) = 0$, $0 < g_* \le g'(u) \le g^* < \infty$, $|g''(u)| \le \bar{g}$ for all $u \in [-\Lambda,\Lambda]$, and $G$ is of the form
\be{ne0}
G = P \circ g,
\ee
where $P$ is a uniformly counterclockwise convex Preisach operator on $[-\Lambda,\Lambda]$. The following result is a straightforward consequence of \cite[Proposition 3.6]{colli}.

\begin{proposition}\label{pc}
	Let $P$ be uniformly counterclockwise convex on $[-\Lambda,\Lambda]$, and let $f$ be an odd increasing function such that $f(0) = 0$. Then there exists $\beta>0$ such that for every sequence $\{w_i: i\in\nat\cup\{-1,0\}\}$ in $[-\Lambda,\Lambda]$ we have
	\be{Pi}
	\begin{aligned}
		&\sumi (P[w]_{i+1} - 2P[w]_i + P[w]_{i-1})f(w_{i+1} - w_i) + \frac{P[w]_0 - P[w]_{-1}}{w_0 - w_{-1}}F(w_0 - w_{-1}) \\
		&\qquad \ge \frac{\beta}{2}\sumi \Gamma(w_{i+1} - w_i),
	\end{aligned}
	\ee
	where we set for $w\in \real$
	\be{gf1}
	F(w) \coloneqq \int_0^w f(v)\dd v, \qquad \Gamma(w) \coloneqq |w|(wf(w) - F(w)) = |w|\int_0^w vf'(v)\dd v.
	\ee
\end{proposition}

We need to define a backward step $u_{-1}$ satisfying the strong formulation of \eqref{dis1} for $i=0$, that is,
\be{e7}
\frac1\tau(G[u]_0(x) - G[u]_{-1}(x)) = \Delta u_0(x) \ \mbox{ in } \Omega
\ee
with homogeneous Neumann boundary condition. Repeating the argument of \cite[Proposition~3.3]{colli} we use Hypotheses \eqref{e3a} and \eqref{c1} to find for each $0<\tau <\rho_0(\Lambda)/2L^2$ functions $u_{-1}$ and $G[u]_{-1}$ satisfying \eqref{e7} as well as the estimate
\be{inim}
\frac1\tau |u_0(x) - u_{-1}(x)| \le C
\ee
with a constant $C>0$ independent of $\tau$ and $x$. The discrete equation \eqref{dis1} extended to $i=0$ has the form
\be{dp1}
\io \left(\frac1\tau(P[w]_i - P[w]_{i-1})\vp + \nabla u_i\cdot\nabla\vp\right)\dd x = 0
\ee
with $w_i = g(u_i)$, for $i\in \nat \cup \{0\}$ and for an arbitrary test function $\vp \in W^{1,2}(\Omega)$. We proceed as in \cite{colli} and test the difference of \eqref{dp1} taken at discrete times $i+1$ and $i$
\be{dp2}
\io \left(\frac1\tau\big(P[w]_{i+1} - 2P[w]_i {+} P[w]_{i-1}\big)\vp + \nabla (u_{i+1}{-} u_{i})\cdot\nabla\vp\right)\dd x = 0
\ee
by $\vp = \tau^{\alpha} f(w_{i+1} - w_i)$ with a suitably chosen odd increasing absolutely continuous function $f$ and exponent $\alpha \in \real$ which will be specified later.

We estimate the initial time increment in the following way.
By \eqref{g} and \eqref{de5a} we have
\be{hatf}
\begin{aligned}
	\tau^{\alpha-1} \frac{P[w]_0 - P[w]_{-1}}{w_0 - w_{-1}}F(w_0 - w_{-1}) &\le \tau^{\alpha-1} \frac{|P[w]_0 - P[w]_{-1}|}{g_*|u_0 - u_{-1}|}F(w_0 - w_{-1}) \\
	&\le C\tau^{\alpha-1} |F(w_0 - w_{-1})| =: \hat F_0(\tau).
\end{aligned}
\ee
Using Proposition~\ref{pc} we have the inequality
\be{dp4}
\tau^{\alpha-1} \sumi (P[w]_{i+1} - 2P[w]_i + P[w]_{i-1})f(w_{i+1} - w_i)
\ge \frac{\beta\tau^{\alpha-1}}{2}\sumi \Gamma(w_{i+1} - w_i) - \hat F_0(\tau)
\ee
with $\hat F_0(\tau)$ defined in \eqref{hatf}. From \eqref{dp2} and \eqref{dp4} we obtain that
\be{dp0}
\frac{\beta\tau^{\alpha-1}}{2}\sumi\io \Gamma (w_{i+1} - w_i)\dd x + \tau^{\alpha}\sumi\io \nabla(u_{i+1}{-}u_i)\cdot \nabla f(w_{i+1}{-} w_i)\dd x \le \hat F_0(\tau).
\ee
Note that
\begin{align*}
	&\nabla(u_{i+1}-u_i)\cdot \nabla f(w_{i+1} - w_i) = f'(w_{i+1} - w_i)\nabla(u_{i+1}-u_i)\cdot \nabla (g(u_{i+1}) - g(u_i))\\
	&\quad = f'(w_{i+1} - w_i)\Big(g'(u_{i+1})|\nabla(u_{i+1}-u_i)|^2 + (g'(u_{i+1}) - g'(u_i))\nabla(u_{i+1}-u_i)\cdot \nabla u_i \Big).
\end{align*}
The properties of $g$ stated in Definition~\ref{dpc} yield
\begin{align} \nonumber
	&\nabla(u_{i+1}-u_i)\cdot \nabla f(w_{i+1} - w_i)\\ \nonumber
	&\quad \ge f'(w_{i+1} - w_i)\Big(g_* |\nabla(u_{i+1}-u_i)|^2 - \bar g |u_{i+1} - u_i|\,|\nabla(u_{i+1}-u_i)|\,|\nabla u_i|\Big)\\ \label{dp6}
	&\quad \ge \frac{g_*}{2}f'(w_{i+1} - w_i)|\nabla(u_{i+1}-u_i)|^2 - Cf'(w_{i+1} - w_i)|u_{i+1}-u_i|^2|\nabla u_i|^2. 
\end{align}
We have
\be{dp7}
g_*|u_{i+1} - u_i| \le |w_{i+1} - w_i|\le g^*|u_{i+1} - u_i|
\ee
for all $i$, and we conclude from \eqref{dp0} that the exists a constant $C>0$ independent of $\tau$ such that
\be{dp8}
\tau^{\alpha-1}\sumi \io \Gamma (w_{i+1} - w_i)\dd x \le C\left(\hat F_0(\tau) + \tau^{\alpha}\sumi\io f'(w_{i+1} - w_i)|w_{i+1}-w_i|^2|\nabla u_i|^2\dd x\right).
\ee
For a suitable $p>1$ we estimate the integral on the right-hand side of \eqref{dp8} using H\"older's inequality as
\begin{align}\nonumber
	&\tau^{\alpha}\sumi\io f'(w_{i+1} - w_i)|w_{i+1}-w_i|^2|\nabla u_i|^2\dd x\\ \label{dp8a}
	&\qquad \le \left(\tau^{1-(1-\alpha) p'}\sumi\io \left(f'(w_{i+1} - w_i)|w_{i+1}-w_i|^2\right)^{p'}\dd x\right)^{1/p'}\left(\tau\sumi\io|\nabla u_i|^{2p} \dd x\right)^{1/p},
\end{align}
where $p' = \frac{p}{p-1}$ is the conjugate exponent to $p$. We now refer to the Gagliardo-Nirenberg inequality in the form
$$
\io|\nabla u_i|^{2p} \dd x \le C\left(\left(\io|\nabla u_i|^{2} \dd x\right)^p + \left(\io|\Delta u_i|^{2} \dd x\right)^{p\kappa}\left(\io|\nabla u_i|^{2} \dd x\right)^{p(1-\kappa)}\right)
$$
with $\kappa = \frac{N}{2p'}$ and a constant $C>0$ independent of $\tau$. We first specify the choice of $p$, namely
\be{kap}
p = 1+\frac{2}{N}.
\ee
Then $p\kappa = 1$, and we obtain
\be{dp9}
\tau\sumi\io|\nabla u_i|^{2p} \dd x \le C\!\left(
\tau\sumi\left(\io|\nabla u_i|^{2} \!\dd x\right)^p\!\! + \tau\sumi\io|\Delta u_i|^{2}\! \dd x \left(\sup_{j\in\nat}\io\!|\nabla u_j|^{2} \dd x\right)^{p-1}\right).
\ee
It follows from \eqref{ea2a} and \eqref{vi} that the right-hand side of \eqref{dp9} is bounded by a constant $C>0$ independent of $\tau$. Using \eqref{dp8}, \eqref{dp8a}, and \eqref{dp9} we thus get the inequality
\begin{align}\nonumber
	&\tau^{\alpha-1}\sumi \io \Gamma (w_{i+1} {-} w_i)\dd x \\ \label{dp10}
	&\qquad\le C\left(\hat F_0(\tau){+} \left(\tau^{1{-}(1-\alpha) p'}\sumi\io \left(f'(w_{i+1} {-} w_i)|w_{i+1}{-}w_i|^2\right)^{p'}\!\!\dd x\right)^{1/p'}\right).
\end{align}
We now consider separately the cases $N \ge 3$ and $N <3$. Assume first that $N\ge 3$. For $p$ from \eqref{kap} and $v \in \real$ we define the functions
\begin{align}\label{fprime}
	f'(v) &= (\tau+|v|)^{p-3},\\ \label{Phi} \Phi(v) &= \int_0^v s f'(s)\dd s = \frac{1}{p-1} \left((\tau + |v|)^{p-1} - \tau^{p-1}\right)- \frac{\tau}{2-p}\left(\tau^{p-2}-(\tau + |v|)^{p-2}\right).
\end{align}
Putting for $z\ge 0$
\begin{align}\label{defq}
	Q(z) &= z \int_0^z s(1+s)^{p-3}\dd s = \frac{z}{p-1}\left((1+z)^{p-1} - 1\right) - \frac{z}{2-p}\left(1 - (1+z)^{p-2}\right),\\[2mm] \label{defm}
	M(z) &= \big(z^2(1+z)^{p-3}\big)^{p'} = z^{2p/(p-1)} (1+z)^{p(p-3)/(p-1)},
\end{align}
and recalling the definition of $\Gamma$ in \eqref{gf1}, we have for $v \in \real$
\be{gamma}
\Gamma(v) = |v|\Phi(v) = \tau^p Q\left(\frac{|v|}{\tau}\right), \quad \left(v^2 f'(v)\right)^{p'} = \tau^p M\left(\frac{|v|}{\tau}\right).
\ee
Choosing $1-\alpha = p-1$, we have $1-(1-\alpha) p' = \alpha-1 = 1-p$. We use \eqref{inim}, \eqref{hatf}, \eqref{dp7}, and the identity
\be{inif}
F(v) = \frac{\tau^{p-1}}{2-p}\left(\frac{|v|}{\tau}- \frac{1}{p-1}\left(\left(1+\frac{|v|}{\tau}\right)^{p-1} - 1\right)\right)
\ee
to check that $\hat F_0(\tau)$ is bounded by a constant $C$ independent of $\tau$, and from \eqref{dp10} and \eqref{gamma} we obtain the estimate
\be{fest}
\tau \sumi \io Q\left(\frac{|w_{i+1} - w_i|}{\tau}\right)\dd x \le C \left(1 + \left(\tau \sumi \io M\left(\frac{|w_{i+1} - w_i|}{\tau}\right)\right)^{1/p'}\dd x\right)
\ee
with a constant $C>0$ independent of $\tau$. We have $2p/(p-1) = N+2$, and formulas \eqref{defq} and \eqref{defm} imply that there exist positive constants $c_1, c_2, c_3, c_4$ such that
\be{asyq}
\lim_{z\to 0} \frac{Q(z)}{z^3} \ge c_1, \quad \lim_{z\to 0} \frac{M(z)}{z^{N+2}} \le c_2, \quad \lim_{z\to \infty} \frac{Q(z)}{z^p} \ge c_3, \quad \lim_{z\to \infty} \frac{M(z)}{z^p} \le c_4. 
\ee
In particular, there exists a constant $K>0$ such that
\be{qm}
M(z) \le K Q(z) \quad \forall z>0.
\ee
The estimate \eqref{fest} thus remains valid if we replace $M$ with $Q$. We conclude that
\be{fq}
\tau \sumi \io Q\left(\frac{|w_{i+1} - w_i|}{\tau}\right)\dd x \le C 
\ee
with a constant $C>0$ independent of $\tau$.

We proceed similarly in the case of dimensions $N<3$. For $N=1$ and $p=3$ we put $f'(v) = 1$ and get the same formulas with
\begin{equation}\label{Q1}
	Q(z) = \frac12 z^3, \quad M(z) = z^3, \quad F(v) = \frac12 |v|^2;
\end{equation}
for $N=2$ and $p=2$ we put $f'(v) = (\tau+|v|)^{-1}$ and a similar computation gives
\begin{equation}\label{Q2}
	Q(z) = z^2 - z \log(1+z), \quad M(z) = z^4/(1+z)^2, \quad F(v) = \tau\Bigg(\bigg(1+\frac{|v|}{\tau}\bigg)\log\!\bigg(1+\frac{|v|}{\tau}\bigg)-\frac{|v|}{\tau}\Bigg),
\end{equation}
with the same conclusion.
In all cases $Q$ is convex with superlinear growth, so that it generates an Orlicz norm on $\Omega$. Estimate \eqref{fq} as the main result of Section \ref{conv} will play a crucial role in the next sections, where we let the discretization parameter $\tau$ tend to $0$ and prove the existence and uniqueness of solutions as well as the asymptotic stabilization result.


\section{Limit as $\tau \to 0$}\label{limi}

We define for $x \in \Omega$ and $t \in [(i-1)\tau, i\tau)$, $i \in \nat$, piecewise linear and piecewise constant interpolations 
\begin{align}
	u\ota(x,t) &= u_{i-1}(x) + \frac{t-(i-1)\tau}{\tau} (u_i(x) - u_{i-1}(x)), \quad \bar{u}\ota(x,t) = u_i(x), \label{inter}\\
	G\ota(x,t) &= G[u]_{i-1}(x) + \frac{t-(i-1)\tau}{\tau} (G[u]_i(x) - G[u]_{i-1}(x)).
\end{align}
Repeating the argument of the proof of \cite[Theorem 1.6]{colli}, we let $\tau$ tend to $0$. On every bounded time interval $(0,T)$ we have bounds independent of $\tau$ for $u\ota_t \in L^p(\Omega\times(0,T))$ by \eqref{fq} and $\nabla u\ota \in L^{2p}(\Omega\times(0,T))$ by \eqref{dp9}. The sequence $u\ota$ is thus compact in $L^p(\Omega; C[0,T])$ and there exists $u\in L^p(\Omega; C[0,T])$ and a subsequence of $u\ota$ such that $u\ota \to u$ and, by Proposition~\ref{pc1}, $G[u\ota] \to G[u]$ and $G\ota \to G[u]$ in $L^p(\Omega; C[0,T])$ strongly for all $T>0$, $u\ota_t \to u_t$ weakly-star in $L^Q(\Omega\times (0,\infty))$ for $Q$ given by \eqref{defq} (or, if $N<3$, by \eqref{Q1} or \eqref{Q2}), $u\ota \to u$ and $\bar{u}\ota \to u$ weakly-star in $L^\infty(0,\infty; W^{1,2}(\Omega))$, and $u$ is the unique solution of \eqref{e0v} satisfying the conditions of Theorem~\ref{t2}. Moreover, using \eqref{Lam}, \eqref{ea2a}, \eqref{vi}, and \eqref{fest} we find constants $\mu>0$, $C>0$ such that
\begin{align}\label{supi}
	\supess \{|u(x,t)|: (x,t) \in \Omega\times(0,\infty)\} &\le \Lambda,\\
	\label{expd}
	\io|\nabla u(x,t)|^2\dd x &\le C\expe^{-\mu t} \quad \mbox{ for all }\ t>0,\\ 
	\label{lapl}
	\int_0^\infty\io |\Delta u(x,t)|^2 \dd x\dd t &\le C,\\
	\label{estime}
	\int_0^\infty\io Q(|u_t(x,t)|)\dd x\dd t &\le C
\end{align}
for $Q$ given by \eqref{defq} with $p = 1+ \frac{2}{N}$. Note that on each bounded time interval $(0,T)$ we can use H\"older's inequality to get from \eqref{asyq} and \eqref{estime} an $L^p$-bound for $u_t$ in the form
\be{lp}
\int_0^T\io |u_t(x,t)|^p \dd x\dd t \le C\left(1 + |\Omega|T\right).
\ee
Finally, still from the fact that \eqref{Lam} is preserved in the limit $\tau \to 0$, that is,
\begin{equation}\label{xi_lim}
	|\xi^r(x,t)| \le (\Lambda - r)^+ \ \mbox{ for }\ale \ x\in\Omega, \mbox{ all } t\ge 0, \mbox{ and all } r\ge 0,
\end{equation}
we obtain that $\xi^r(x,t) = 0$ for all $r\ge \Lambda$. Hence, similarly as in \eqref{Gi_Lambda}, we actually have
\begin{equation}\label{G_Lambda}
	G[u](x,t) = \bar G + \ \int_0^\Lambda \psi(x,r,\xi^r(x,t)) \dd r.
\end{equation}


\section{Asymptotics as $t \to \infty$}\label{sec:asym}

Notice first that choosing in relation \eqref{e4a} $z = r\sign(\xi^r_t)$, we get $r|\xi^r_t|\le \xi^r_t (u - \xi^r) \le r|\xi^r_t|$ a.\,e., hence, $\xi^r_t (u - \xi^r) = r|\xi^r_t|$ a.\,e.
With the notation of \eqref{psi}, we have
\begin{align*}
	\rho(x,r,\xi^r) r |\xi^r_t| &= \rho(x,r,\xi^r) \xi^r_t (u - \xi^r) =  (u - \xi^r)\frac{\partial}{\partial t}\psi(x,r,\xi^r)\\ & = u \frac{\partial}{\partial t}\psi(x,r,\xi^r)- \frac{\partial}{\partial t} \big(\psi(x,r,\xi_r)\xi_r - \Psi(x,r,\xi_r)\big) \ \ale
\end{align*}
which gives the hysteresis energy balance equation
\be{ene}
G[u]_t u = \frac{\partial}{\partial t}\int_0^\Lambda \big(\psi(x,r,\xi^r)\xi^r - \Psi(x,r,\xi^r)\big) \dd r + \int_0^\Lambda \rho(x,r,\xi^r) r|\xi^r_t| \dd r,
\ee
with $G[u]$ as in $\eqref{G_Lambda}$. To prove the convergence of $\xi^r$ (and, consequently, the convergence of $G[u]$) as $t\to \infty$, we test \eqref{e0v} by $\vp = u$ and
use \eqref{ene} to obtain (omitting the arguments $x$ and $r$ of $\rho$, $\psi$, and $\Psi$ for simplicity)
\be{ener}
\frac{\dd}{\dd t} \io\int_0^\Lambda \big(\psi(\xi^r)\xi^r - \Psi(\xi^r)\big)\dd r\dd x
+ \io\int_0^\Lambda \rho(\xi^r)r|\xi^r_t|\dd r\dd x +\io |\nabla u|^2\dd x = 0.
\ee
By \eqref{e3a} we have $\psi(\xi)\xi - \Psi(\xi) = \int_0^\xi s\rho(s)\dd s \ge \frac12\rho_0(\Lambda)|\xi|^2 \ge 0$ for every $|\xi|\le \Lambda$, which is our case thanks to \eqref{xi_lim}. Hence, integrating \eqref{ener} over $t \in (0,\infty)$ we get
\be{varxi}
\int_0^\infty\io\int_0^\Lambda r|\xi^r_t|\dd r\dd x\dd t \le C
\ee
with a constant $C>0$ depending only on the initial condition. For every sequence $0 = t_0 < t_1 < t_2 < \dots$ we thus have
\be{varj}
\sum_{j=1}^\infty\io\int_0^\Lambda r|\xi^r(x,t_j) - \xi^r(x,t_{j-1}) |\dd r\dd x \le C.
\ee
The sequence $\{\xi^r(x,t_j)\}$ is therefore a fundamental (Cauchy) sequence in the space $L^1(\Omega\times(0,\Lambda))$ endowed with the weighted norm $\|\zeta\| = \io\int_0^\Lambda r|\zeta(x,r)|\dd x\dd r$. This is indeed a Banach space, and we conclude that there exists $\bar\lambda \in L^\infty(\Omega\times (0,\Lambda))$, $\bar\lambda(x,r) = 0$ for $r>\Lambda$, $|\bar\lambda(x,r) - \bar\lambda(x,s)| \le r-s$ for all $r>s\ge 0$ such that 
\be{asy}
\bar\omega(t) \coloneqq \io\int_0^\Lambda r|\xi^r(x,t) - \bar\lambda(x,r)|\dd r\dd x \to 0 \quad \mbox {as } \ t\to \infty.
\ee
Thanks to the above convergence, we can choose $T>0$ such that $\bar\omega(t)\le (\Lambda/3)^3$ for $t\ge T$. For $r\in [0,\Lambda]$ and $t\ge T$ put
$$
\omega(r,t) \coloneqq \io |\xi^r(x,t) - \bar\lambda(x,r)|\dd x,
$$
with the intention to prove that $\omega(r,t) \to 0$ as $t\to\infty$ uniformly with respect to $r$, and that the convergence rate can be estimated in terms of $\bar{\omega}(t)$ given by \eqref{asy}. To this aim, for $t \ge T$ we introduce the sets
$$
A(t) \coloneqq \{r\ge \bar\omega(t)^{1/3}: \ \omega(r,t) \ge \bar\omega(t)^{1/3}\}.
$$
Using \eqref{asy} we get
$$
\bar\omega(t)  = \int_0^\Lambda r \omega(r,t)\dd r \ge \int_{A(t)} r \omega(r,t)\dd r \ge \bar\omega(t)^{2/3} |A(t)|,
$$
which yields the upper bound for the Lebesgue measure of $A(t)$ in the form
\begin{equation}\label{misA}
	|A(t)| \le \bar\omega(t)^{1/3}.
\end{equation}
Let now $r\in (0,\Lambda)$ be arbitrary. In the ``good'' case $r \in [\bar\omega(t)^{1/3},\Lambda] \setminus A(t)$ we immediately have by definition  of $A(t)$ the desired bound $\omega(r,t) < \bar\omega(t)^{1/3}$. Instead, in the ``bad'' cases $r\in (0,\bar\omega(t)^{1/3})$ or $r \in A(t)$, by \eqref{misA} we find a ``good'' $s\in [\bar\omega(t)^{1/3},\Lambda]\setminus A(t)$ such that $|r-s| \le 2\bar\omega(t)^{1/3}$, and thanks to \eqref{lip-r} we estimate
$$
\omega(r,t) \le \omega(s,t) + 2|\Omega||r-s| \le (1+4|\Omega|)\bar\omega(t)^{1/3}\ \for t\ge T.
$$
This shows that $\xi^r(x,t) \to \bar\lambda(x,r)$ strongly in $L^1(\Omega)$ and uniformly with respect to $r$ as $t\to\infty$. By the Lebesgue dominated convergence theorem we also obtain the strong convergence in $L^q(\Omega)$ for all $q \ge 1$, still uniformly in $r$.

To prove the convergence of $u$, we define the mean value of $u(x,t)$
\be{mean}
U(t) \coloneqq \frac{1}{|\Omega|}\io u(x,t)\dd x.
\ee
For all $t>0$ we have $|U(t)| \le \Lambda$ by virtue of \eqref{supi}, and 
\be{pw} 
\io |u(x,t) - U(t)|^2\dd x \le C \io|\nabla u(x,t)|^2\dd x
\ee
with some constant $C>0$ by virtue of the classical Poincar\'e-Wirtinger inequality. From \eqref{expd} it follows that
\be{pw0} 
\io |u(x,t) - U(t)|^2\dd x \le C \expe^{-\mu t}
\ee
with some constants $\mu>0$ and $C>0$. To conclude the proof of Theorem~\ref{t2} it suffices to check that there exists a constant $\bar u \in \real$ such that
\be{Ua}
\lim_{t \to \infty} U(t) = \bar u.
\ee
The fact that $\bar\lambda(x,0) = \bar u$ then follows from the uniform (with respect to $r$) convergence $\xi^r(x,t) \to \bar\lambda(x,r)$.

To prove \eqref{Ua}, we proceed by contradiction. Assume that \eqref{Ua} does not hold. In this case, we would have
\be{Uac}
\liminf_{t \to \infty} U(t) = A < B = \limsup_{t \to \infty} U(t).
\ee
Using \eqref{pw0} we find an increasing sequence $\{t_j\}$, $t_j \to \infty$ as $j \to \infty$, such that
\begin{align}
	&\io |u(x,t_j) - U(t_j)|^2\dd x \le \expe^{-2j} \ \ \forall j\in\nat, \label{pw1}\\
	&U(t_{2i}) \to A, \ U(t_{2i+1}) \to B \ \mbox{ as } i\to \infty. \label{pw2}
\end{align}
We define the sets $\Omega_j = \{x \in \Omega: |u(x,t_j) - U(t_j)|^2 > \expe^{-j}\}$. By \eqref{pw1}, we have $|\Omega_j| \le \expe^{-j}$. Consider now $j\ge n$ for $n\in\nat$, and put
$$
\Omega\on \coloneqq \bigcup_{j=n}^{\infty} \Omega_j. 
$$
Then $|\Omega\on| \le \frac{\expe}{\expe-1} \expe^{-n}$, and by \eqref{pw2} we may assume, choosing $n$ sufficiently large, that for some $A < a < b < B$, for a.\,e.\ $x \in \Omega\setminus \Omega\on$, and for $i = 1,2,...$ we have
\be{ii}
u(x, t_{2i}) \le a, \quad u(x,t_{2i+1}) \ge b.
\ee
From the elementary inequalities
$$
\xi^r(x,t_{2i}) \le u(x, t_{2i})+ r, \quad \xi^r(x,t_{2i+1}) \ge u(x, t_{2i+1}) - r
$$
and from \eqref{ii} we infer that
\be{ii1}
\xi^r(x,t_{2i+1})-\xi^r(x,t_{2i}) \ge b-a-2r > 0 
\ee
for $r\in (0,(b-a)/2)$ and $x \in \Omega\setminus \Omega\on$. Hence, for $r < (b-a)/2$ we get
$$
\io |\xi^r(x,t_{2i+1}){-}\xi^r(x,t_{2i})|\dd x \ge \int_{\Omega\setminus \Omega\on}\!|\xi^r(x,t_{2i+1}){-}\xi^r(x,t_{2i})|\dd x \ge (|\Omega| {-} |\Omega\on|) (b{-}a{-}2r).
$$
We thus have
\begin{align*}
	\io\int_0^\Lambda r |\xi^r(x,t_{2i+1})-\xi^r(x,t_{2i})|\dd x \dd r
	&\ge (|\Omega| - |\Omega\on|) \int_0^{\min\{\Lambda,(b-a)/2\}} r(b-a-2r)\dd r \\
	&\ge (|\Omega| - |\Omega\on|)\left(\min\{\Lambda,(b-a)/2\}\right)^2\,\frac{b-a}{6} > 0
\end{align*}
for all $i \ge n$, which contradicts \eqref{varj}. We thus conclude that \eqref{Ua} holds, which completes the proof of Theorem~\ref{t2}.


\end{document}